\documentclass[12pt,reqno]{amsart}
\usepackage{amssymb}
\usepackage{epsfig}
\usepackage{graphicx}

\newtheorem{definition}{Definition}[section]
\newtheorem{theorem}{Theorem}[section]
\newtheorem{lemma}{Lemma}[section]

\newtheorem{prop}{Proposition}[section]
\newtheorem{remark}{Remark}[section]

\newcommand{\bt}{\begin{theorem}}
\newcommand{\et}{\end{theorem}}
\newcommand{\bl}{\begin{lemma}}
\newcommand{\el}{\end{lemma}}
\newcommand{\bd}{\begin{definition}}
\newcommand{\ed}{\end{definition}}
\newcommand{\bc}{\begin{corollary}}
\newcommand{\ec}{\end{corollary}}
\newcommand{\bp}{\begin{proof}}
\newcommand{\ep}{\end{proof}}
\newcommand{\bx}{\begin{example}}
\newcommand{\ex}{\end{example}}
\newcommand{\bi}{\begin{exercise}}
\newcommand{\ei}{\end{exercise}}
\newcommand{\bo}{\begin{prop}}
\newcommand{\eo}{\end{prop}}
\newcommand{\br}{\begin{remark}}
\newcommand{\er}{\end{remark}}
\newcommand{\be}{\begin{equation}}
\newcommand{\ee}{\end{equation}}
\newcommand{\ba}{\begin{align}}
\newcommand{\ea}{\end{align}}
\newcommand{\bn}{\begin{enumerate}}
\newcommand{\en}{\end{enumerate}}
\newcommand{\bg}{\begin{align*}}
\newcommand{\bcs}{\begin{cases}}
\newcommand{\ecs}{\end{cases}}

\newcommand{\bean}{\begin{eqnarray*}}
\newcommand{\eean}{\end{eqnarray*}}

\numberwithin{equation}{section}
\begin{document}
\title{{\bf Nondegeneracy of half-harmonic maps from $\mathbb{R}$ into $\mathbb{S}^1$}}

\author[Y. Sire]{Yannick Sire}
\address{Johns Hopkins University, Department of mathematics, Krieger Hall,  Baltimore, MD 21218, USA}
\email{sire@math.jhu.edu}

\author[J. Wei]{Juncheng Wei}
\address{Department of
Mathematics,  University of British Columbia, Vancouver, B.C.,
Canada, V6T 1Z2}
\email{wei@math.cuhk.edu.hk}

\author[Y.Zheng]{Youquan ZHENG}
\address{School of Science,
Tianjin University,
92 Weijin Road, Tianjin 300072,
P.R. China}
\email{zhengyq@tju.edu.cn}

\date{}
\maketitle
\begin{center}
\begin{minipage}{120mm}
\begin{center}{\bf Abstract}\end{center}
We prove that the standard half-harmonic map $U:\mathbb{R}\to\mathbb{S}^1$ defined by
\begin{equation*}
x\to \begin{pmatrix}
     \frac{x^2-1}{x^2+1} \\
     \frac{-2x}{x^2+1}
\end{pmatrix}
\end{equation*}
is nondegenerate in the sense that all bounded solutions of the linearized half-harmonic map equation are linear combinations of
three functions corresponding to rigid motions (dilation, translation and rotation) of $U$.

\end{minipage}
\end{center}
\vskip0.10in
\section{Introduction}
Due to their importance in geometry and physics, the analysis of critical points of conformal invariant Lagrangians has attracted
much attention since 1950s. A typical example is the Dirichlet energy which is defined on two-dimensional domains and its critical points are harmonic maps. This definition can be generalized to even-dimensional domains whose critical points are called polyharmonic maps. In recent years, people are very interested in the analog of Dirichlet energy in odd-dimensional case, for example, \cite{DaLioCVPDE2013},
\cite{DaLioAIHAN2015}, \cite{DaLioAdvMath2011}, \cite{DaLioAnalPDE2011}, \cite{MillotSireARMA2015}, \cite{SchikorraJDE2012} and the references therein. Among these works, a special case is the so-called half-harmonic maps from $\mathbb{R}$ into $\mathbb{S}^{1}$ which are defined as critical points of the line energy
\begin{equation}\label{e:linearenergy}
\mathcal{L}(u) = \frac{1}{2}\int_{\mathbb{R}}|(-\Delta_{\mathbb{R}})^{\frac{1}{4}}u|^2dx.
\end{equation}
Note that the functional $\mathcal{L}$ is invariant under the
trace of conformal maps keeping invariant the half-space $\mathbb{R}^2_+$: the M\"{o}bius group. Half-harmonic maps have close relations with harmonic maps with partially free boundary and minimal surfaces with free boundary, see \cite{FraserSchoenInventMath2016} and \cite{MillotSireARMA2015}. Computing the associated Euler-Lagrange equation of (\ref{e:linearenergy}), we obtain that if $u:\mathbb{R}\to \mathbb{S}^1$ is a half-harmonic map, then $u$ satisfies the following equation,
\begin{equation}\label{e:main}
(-\Delta_{\mathbb{R}})^{\frac{1}{2}}u(x) = \left(\frac{1}{2\pi}\int_{\mathbb{R}}\frac{|u(x)-u(y)|^2}{|x-y|^2}dy\right)u(x)\text{ in }\mathbb{R}.
\end{equation}
It was proved in \cite{MillotSireARMA2015} that
\begin{prop}(\cite{MillotSireARMA2015})
Let $u\in \dot{H}^{1/2}(\mathbb{R},\mathbb{S}^1)$ be a non-constant entire half-harmonic map into $\mathbb{S}^1$ and $u^e$ be its harmonic extension to $\mathbb{R}^2_+$. Then there exist $d\in\mathbb{N}$, $\vartheta\in \mathbb{R}$, $\{\lambda_k\}_{k=1}^d\subset (0, \infty)$ and $\{a_k\}_{k=1}^d\subset \mathbb{R}$ such that $u^e(z)$ or its complex conjugate equals to
$$
e^{i\vartheta}\prod_{k=1}^d\frac{\lambda_k(z-a_k)-i}{\lambda_k(z-a_k)+i}.
$$
Furthermore,
$$
\mathcal{E}(u,\mathbb{R}) = [u]^2_{H^{1/2}(\mathbb{R})}=\frac{1}{2}\int_{\mathbb{R}^2_+}|\nabla u^e|^2dz = \pi d.
$$
\end{prop}
This proposition shows that the map $U:\mathbb{R}\to\mathbb{S}^1$
\begin{equation*}
x\to \begin{pmatrix}
     \frac{x^2-1}{x^2+1} \\
     \frac{-2x}{x^2+1}
\end{pmatrix}
\end{equation*}
is a half-harmonic map corresponding to the case $\vartheta = 0$, $d=1$, $\lambda_1=1$ and $a_1=0$.
In this paper, we prove the nondegeneracy of $U$ which is a crucial ingredient when analyzing the singularity formation of half-harmonic map flow.
Note that $U$ is invariant under translation, dilation and rotation, i.e., for $Q=\begin{pmatrix}
     \cos\alpha & -\sin\alpha \\
     \sin\alpha & \cos\alpha
\end{pmatrix}\in O(2)$, $q\in \mathbb{R}$ and $\lambda\in\mathbb{R}^+$, the function
\begin{equation*}
QU\left(\frac{x-q}{\lambda}\right) = \begin{pmatrix}
     \cos\alpha & -\sin\alpha \\
     \sin\alpha & \cos\alpha
\end{pmatrix}U\left(\frac{x-q}{\lambda}\right)
\end{equation*}
still satisfies (\ref{e:main}). Differentiating with $\alpha$, $q$ and $\lambda$ respectively and then set $\alpha = 0$, $q=0$ and $\lambda = 1$, we obtain that the following three functions
\begin{equation}\label{e:kernel}
Z_1(x) =
\begin{pmatrix}
     \frac{2x}{x^2+1} \\
     \frac{x^2-1}{x^2+1}
\end{pmatrix},\quad
Z_2(x) = \begin{pmatrix}
     \frac{-4x}{(x^2+1)^2} \\
     \frac{2(1-x^2)}{(x^2+1)^2}
\end{pmatrix},\quad
Z_3(x) =
\begin{pmatrix}
     \frac{-4x^2}{(x^2+1)^2} \\
     \frac{2x(1-x^2)}{(x^2+1)^2}
\end{pmatrix}
\end{equation}
satisfy the linearized equation at the solution $U$ of (\ref{e:main}) defined as
\begin{eqnarray}\label{e:linearized1}
\nonumber (-\Delta_{\mathbb{R}})^{\frac{1}{2}}v(x) &=& \left(\frac{1}{2\pi}\int_{\mathbb{R}}\frac{|U(x)-U(y)|^2}{|x-y|^2}dy\right)v(x)\\
\quad\quad\quad &&+ \left(\frac{1}{\pi}\int_{\mathbb{R}}\frac{(U(x)-U(y))\cdot(v(x) -v(y))}{|x-y|^2}dy\right)U(x)\quad\text{in }\mathbb{R}
\end{eqnarray}
for $v:\mathbb{R}\to T_U\mathbb{S}^1$.
Our main result is
\begin{theorem}\label{t:main}
The half-harmonic map $U:\mathbb{R}\to\mathbb{S}^1$
\begin{equation*}
x\to \begin{pmatrix}
     \frac{x^2-1}{x^2+1} \\
     \frac{-2x}{x^2+1}
\end{pmatrix}
\end{equation*}
is nondegenerate in the sense that all bounded solutions of equation (\ref{e:linearized1}) are linear combinations of $Z_1$, $Z_2$ and $Z_3$ defined in (\ref{e:kernel}).
\end{theorem}

In the case of harmonic maps from two-dimensional domains into $\mathbb{S}^2$, the non-degeneracy of bubbles was proved in Lemma 3.1 of \cite{DelPinoWeiPreprint}.
Integro-differential equations have attracted substantial research in recent
years. The nondegeneracy of ground state solutions for the fractional nonlinear Schr\"{o}dinger equations has been proved by Frank and
Lenzmann \cite{FrankLenzmannActaMath2013}, Frank, Lenzmann and Silvestre \cite{FrankLenzmannCPAM}, Fall and Valdinoci \cite{FallCMP2014}, and the corresponding result in the case of fractional Yamabe problem was obtained by D{\'a}vila, del Pino and Sire in \cite{DelPinoSirePAMS}.

\section{Proof of Theorem \ref{t:main}}
The rest of this paper is devoted to the proof of Theorem \ref{t:main}. For convenience, we identify $\mathbb{S}^1$ with the complex unite circle. Since $Z_1$, $Z_2$ and $Z_3$ are linearly independent and belong to the space $L^\infty(\mathbb{R})\cap Ker(\mathcal{L}_0)$, we only need to prove that the dimension of $L^\infty(\mathbb{R})\cap Ker(\mathcal{L}_0)$ is 3. Here the operator $\mathcal{L}_0$ is defined as
\begin{eqnarray*}
&&\mathcal{L}_0(v)=(-\Delta_{\mathbb{R}})^{\frac{1}{2}}v(x) - \left(\frac{1}{2\pi}\int_{\mathbb{R}}\frac{|U(x)-U(y)|^2}{|x-y|^2}dy\right)v(x)\\
&&\nonumber\quad\quad\quad\quad - \left(\frac{1}{\pi}\int_{\mathbb{R}}\frac{(U(x)-U(y))\cdot(v(x) -v(y))}{|x-y|^2}dy\right)U(x),
\end{eqnarray*}
for $v:\mathbb{R}\to T_U\mathbb{S}^1$.
Let us come back to equation (\ref{e:linearized1}), for $v:\mathbb{R}\to T_U\mathbb{S}^1$, $v(x)\cdot U(x) = 0$ holds pointwisely. Using this fact and the definition of $(-\Delta_{\mathbb{R}})^{\frac{1}{2}}$ (see \cite{NPVBSM2012}), we have
\begin{eqnarray*}
(-\Delta_{\mathbb{R}})^{\frac{1}{2}}v(x) &=& \left(\frac{1}{2\pi}\int_{\mathbb{R}}\frac{|U(x)-U(y)|^2}{|x-y|^2}dy\right)v(x)\\
\quad\quad\quad &&+ \left(\frac{1}{\pi}\int_{\mathbb{R}}\frac{(U(x)-U(y))\cdot(v(x) -v(y))}{|x-y|^2}dy\right)U(x)\\
&= &\left(\frac{1}{2\pi}\int_{\mathbb{R}}\frac{|U(x)-U(y)|^2}{|x-y|^2}dy\right)v(x)\\
\quad\quad\quad &&+ \left(\frac{1}{\pi}\int_{\mathbb{R}}\frac{(U(x)-U(y))}{|x-y|^2}dy\cdot v(x)\right)U(x)\\
\quad\quad\quad &&+ \left(\frac{1}{\pi}\int_{\mathbb{R}}\frac{(v(x)-v(y))}{|x-y|^2}dy\cdot U(x)\right)U(x)\\
&= &\left(\frac{1}{2\pi}\int_{\mathbb{R}}\frac{|U(x)-U(y)|^2}{|x-y|^2}dy\right)v(x)\\
\quad\quad\quad &&+ \left(\frac{1}{\pi}\int_{\mathbb{R}}\frac{(v(x)-v(y))}{|x-y|^2}dy\cdot U(x)\right)U(x)\\
&= &\left(\frac{1}{2\pi}\int_{\mathbb{R}}\frac{|U(x)-U(y)|^2}{|x-y|^2}dy\right)v(x)\\
\quad\quad\quad &&+ \left((-\Delta_{\mathbb{R}})^{\frac{1}{2}}v(x)\cdot U(x)\right)U(x).
\end{eqnarray*}
Therefore equation (\ref{e:linearized1}) becomes to
\begin{eqnarray}\label{e:aftertransform}
\nonumber (-\Delta_{\mathbb{R}})^{\frac{1}{2}}v(x) &=& \left(\frac{1}{2\pi}\int_{\mathbb{R}}\frac{|U(x)-U(y)|^2}{|x-y|^2}dy\right)v(x) + \left((-\Delta_{\mathbb{R}})^{\frac{1}{2}}v(x)\cdot U(x)\right)U(x)\\
&=&\frac{2}{x^2+1}v(x) + \left((-\Delta_{\mathbb{R}})^{\frac{1}{2}}v(x)\cdot U(x)\right)U(x).
\end{eqnarray}

Next, we will lift equation (\ref{e:aftertransform}) to $\mathbb{S}^1$ via the stereographic projection from $\mathbb{R}$ to $\mathbb{S}^1\setminus\{pole\}$:
\begin{equation}\label{e:stereographicprojection}
S(x) = \begin{pmatrix}
     \frac{2x}{x^2+1} \\
     \frac{1-x^2}{x^2+1}
\end{pmatrix}.
\end{equation}
It is well known that  the Jacobian of the stereographic projection is
$$
J(x) =  \frac{2}{x^2+1}.
$$
For a function $\varphi:\mathbb{R}\to \mathbb{R}$, define $\tilde{\varphi}:\mathbb{S}^1\to \mathbb{R}$ by
\begin{equation}\label{e:lift}
\varphi(x) = J(x)\tilde{\varphi}(S(x)).
\end{equation}
Then we have
\begin{eqnarray*}
[(-\Delta_{\mathbb{S}^1})^{\frac{1}{2}}\tilde{\varphi}](S(x)) &=& \frac{1}{\pi}\int_{\mathbb{R}}\frac{\tilde{\varphi}(S(x))- \tilde{\varphi}(S(y))}{|S(x)-S(y)|^2}dS(y)\\
&=& \frac{1}{\pi}\int_{\mathbb{R}}\frac{\frac{1+x^2}{2}\varphi(x)- \frac{1+y^2}{2}\varphi(y)}{\frac{4(x-y)^2}{(x^2+1)(y^2+1)}}\frac{2}{1+y^2}dy\\
&=& \frac{1+x^2}{4\pi}\int_{\mathbb{R}}\frac{(1+x^2)\varphi(x)- (1+y^2)\varphi(y)}{(x-y)^2}dy\\
&=& \frac{1+x^2}{2}(-\Delta_{\mathbb{R}})^{1/2}\left[\frac{x^2+1}{2}\varphi(x)\right]\\
&=& \frac{1+x^2}{2}(-\Delta_{\mathbb{R}})^{1/2}\left[\tilde{\varphi}(S(x))\right].
\end{eqnarray*}
Therefore,
$$
(-\Delta_{\mathbb{R}})^{1/2}\left[\tilde{\varphi}(S(x))\right] = J(x)[(-\Delta_{\mathbb{S}^1})^{\frac{1}{2}}\tilde{\varphi}](S(x)).
$$

Denote $v = (v_1, v_2)$ and let $\tilde{v}_1$, $\tilde{v}_2$ be the functions defined by (\ref{e:lift}) respectively. Then the linearized equation (\ref{e:aftertransform}) becomes
\begin{equation*}
\left\{\begin{array}{ll}
        J(x)(-\Delta_{\mathbb{S}^1})^{\frac{1}{2}}\tilde{v}_1 = J(x)\tilde{v}_1 + \frac{x^2-1}{x^2+1}\frac{x^2-1}{x^2+1}J(x)(-\Delta_{\mathbb{S}^1})^{\frac{1}{2}}\tilde{v}_1 + \frac{x^2-1}{x^2+1}\frac{-2x}{x^2+1}J(x)(-\Delta_{\mathbb{S}^1})^{\frac{1}{2}}\tilde{v}_2,\\
        J(x)(-\Delta_{\mathbb{S}^1})^{\frac{1}{2}}\tilde{v}_2 = J(x)\tilde{v}_2 + \frac{-2x}{x^2+1}\frac{x^2-1}{x^2+1}J(x)(-\Delta_{\mathbb{S}^1})^{\frac{1}{2}}\tilde{v}_1 + \frac{-2x}{x^2+1}\frac{-2x}{x^2+1}J(x)(-\Delta_{\mathbb{S}^1})^{\frac{1}{2}}\tilde{v}_2.
       \end{array}
\right.
\end{equation*}
Since $J(x) > 0$ and set $U = (\cos\theta, \sin\theta)$, we get
\begin{equation*}
\left\{\begin{array}{ll}
        (-\Delta_{\mathbb{S}^1})^{\frac{1}{2}}\tilde{v}_1 = \tilde{v}_1 + \cos^2\theta(-\Delta_{\mathbb{S}^1})^{\frac{1}{2}}\tilde{v}_1 + \cos\theta\sin\theta(-\Delta_{\mathbb{S}^1})^{\frac{1}{2}}\tilde{v}_2,\\
        (-\Delta_{\mathbb{S}^1})^{\frac{1}{2}}\tilde{v}_2 = \tilde{v}_2 + \cos\theta\sin\theta(-\Delta_{\mathbb{S}^1})^{\frac{1}{2}}\tilde{v}_1 +\sin^2\theta(-\Delta_{\mathbb{S}^1})^{\frac{1}{2}}\tilde{v}_2,
       \end{array}
\right.
\end{equation*}
which is equivalent to
\begin{equation*}
\left\{\begin{array}{ll}
        (-\Delta_{\mathbb{S}^1})^{\frac{1}{2}}\tilde{v}_1 = 2\tilde{v}_1 + \cos2\theta(-\Delta_{\mathbb{S}^1})^{\frac{1}{2}}\tilde{v}_1 + \sin2\theta(-\Delta_{\mathbb{S}^1})^{\frac{1}{2}}\tilde{v}_2,\\
        (-\Delta_{\mathbb{S}^1})^{\frac{1}{2}}\tilde{v}_2 = 2\tilde{v}_2 + \sin2\theta(-\Delta_{\mathbb{S}^1})^{\frac{1}{2}}\tilde{v}_1 -\cos2\theta(-\Delta_{\mathbb{S}^1})^{\frac{1}{2}}\tilde{v}_2.
       \end{array}
\right.
\end{equation*}
Set $w = \tilde{v}_1 + i\tilde{v}_2$, $z=\cos\theta + i\sin\theta$, then we have
\begin{equation}\label{e:linearizedfinal}
(-\Delta_{\mathbb{S}^1})^{\frac{1}{2}}w = 2w + z^2(-\Delta_{\mathbb{S}^1})^{\frac{1}{2}}\bar{w}.
\end{equation}
Here $\bar{w}$ is the conjugate of $w$.

Since $v\in L^\infty(\mathbb{R})$, $w$ is also bounded, so we can expand $w$ into fourier series
$$
w = \sum_{k=-\infty}^\infty a_kz^k.
$$
Note that all the eigenvalues for $(-\Delta_{\mathbb{S}^1})^{\frac{1}{2}}$ are $\lambda_k = k$, $k = 0, 1, 2, \cdot\cdot\cdot$, see \cite{ChangYangAnnals1995}. Using (\ref{e:linearizedfinal}), $(-\Delta_{\mathbb{S}^1})^{\frac{1}{2}}z^k = kz^k$ and $(-\Delta_{\mathbb{S}^1})^{\frac{1}{2}}\bar{z}^k = k\bar{z}^k$, we obtain
\begin{equation*}
\left\{\begin{array}{ll}
        (-k-2)a_k = (2-k)\bar{a}_{2-k}, \text{ if } k<0,\\
        (k-2)a_k = (2-k)\bar{a}_{2-k}, \text{ if } 0\leq k\leq 2,\\
        a_k=\bar{a}_{2-k}, \text{ if } k\geq 3.
       \end{array}
\right.
\end{equation*}
Furthermore, from the orthogonal condition $v(x)\cdot U(x) = 0$ (so $(\tilde{v}_1, \tilde{v}_2)\cdot (\cos\theta, \sin\theta) = 0$), we have
$$
a_k = -\bar{a}_{2-k}, \quad k=\cdots -1, 0, 1, \cdots.
$$
Thus
$$
a_k = 0, \text{ if } k<0 \text{ or } k\geq 3
$$
and
$$
a_0 = -\bar{a}_2, \quad a_1=-\bar{a}_1
$$
hold, which imply that
$$
w = -\bar{a}_2 + a_1z + a_2z^2 = a(iz) + b\left[\frac{i}{2}(z-1)^2\right] + c\frac{(z^2-1)}{2}.
$$
Here $a$, $b$, $c$ are real numbers and satisfy relations
$$
i(a-b)=a_1,\quad \frac{c}{2} + \frac{i}{2}b = a_2.
$$
And it is easy to check that $iz$, $\frac{i}{2}(z-1)^2$ and $\frac{(z^2-1)}{2}$ are respectively $Z_1$, $Z_2$ and $Z_3$ under stereographic projection (\ref{e:stereographicprojection}). By the one-to-one correspondence of $w$ and $v$, we know that the dimension of $L^\infty(\mathbb{R})\cap Ker(\mathcal{L}_0)$ is 3.
This completes the proof.

\end{document}